\documentclass[journal]{IEEEtran}
\newcommand{\subparagraph}{}
\usepackage{titlesec}
\usepackage{textcomp}
\usepackage{eurosym}
\usepackage[english]{babel}
\usepackage{float}
\usepackage[utf8x]{inputenc}
\usepackage{bm}
\usepackage{enumitem}
\usepackage{booktabs}
\usepackage[table,xcdraw]{xcolor}

\usepackage{amsmath}
\usepackage{amsfonts}
\usepackage{comment}
\usepackage[lined]{algorithm2e}
\usepackage{amsthm}
\usepackage{graphicx}
\usepackage{lettrine}
\usepackage[colorlinks=true, allcolors=black]{hyperref}

\newtheorem{thm}{Theorem}
\newtheorem{lemma}{Lemma}

\theoremstyle{definition}

\newcommand*{\QEDB}{\hfill\ensuremath{\square}}%
\newcommand{\bd}{\boldsymbol}

\makeatletter
\renewcommand\paragraph{\@startsection{paragraph}{4}{\z@}%
    {1ex \@plus0ex \@minus0ex}%
    {-1em}%
    {\normalfont\normalsize\bfseries}}
\makeatother

\DeclareRobustCommand{\officialeuro}{%
	\ifmmode\expandafter\text\fi
	{\fontencoding{U}\fontfamily{eurosym}\selectfont e}}

\usepackage{footnote}
\makesavenoteenv{tabular}
\makesavenoteenv{table}
\ifCLASSINFOpdf

\fi

\begin{document}

\title{Non-Wire Alternatives: an Additional Value Stream for Distributed Energy Resources }

\author{\IEEEauthorblockN{Jesus E. Contreras-Oca\~na, Yize Chen, Uzma Siddiqi, and Baosen Zhang }
\thanks{J. Contreras-Oca\~na is with the Universit\'e Grenoble Alpes  and the Laboratoire de G\'enie El\'ectrique de Grenoble (G2Elab), (email: jesus.contreras.ocana@gmail.com); Y. Chen and B. Zhang are with ECE at the University of Washington (emails: \{yizechen, zhangbao\}@uw.edu); U. Siddiqi is with Seattle City Light (email:Uzma.Siddiqi@seattle.gov).}
}

\maketitle

\begin{abstract}
    Distributed energy resources (DERs) can serve as non-wire alternatives (NWAs) to capacity expansion by managing peak load to avoid or delay traditional expansion projects. However, the value stream derived from using DERs as NWAs is usually not explicitly included in DER planning problems. In this paper, we study a planning problem that co-optimizes investment and operation of DERs and the timing of capacity expansion. By including the timing of capacity expansion as a decision variable, we naturally incorporate NWA value stream of DERs into the planning problem. Furthermore,  we show that even though the resulting optimization problem could have millions of variables and is non-convex, an optimal solution can be found by solving a series of smaller linear problems.  Finally, we present a NWAs planning problem using real data from the Seattle Campus of the University of Washington. 
\end{abstract}

\IEEEpeerreviewmaketitle

\section{Introduction}
Power systems are typically designed for the peak load, which happens a small number of hours per year. When the load reaches capacity, the traditional solution is to expand generation capacity, install more wires, or reinforce existing ones~\cite{seifi2011electric, GACITUA2018346}. While decades of experience makes this conventional or ``wires'' solution reliable and safe, it often carries high capital costs, can face public opposition, and experience time-consuming legal issues (e.g., eminent domain questions)~\cite{stanton2015getting}.

Lately, regulators have pushed system planners to consider distributed energy resources (DERs), network management, grid optimization, dynamic pricing, and data analytics as alternative methods to the “wires” solution.  For example, the I-5 corridor project in the Pacific Northwest of the United States~\cite{John2017_NWA_BPA} explores alternatives to transmission capacity and the Brooklyn-Queens Demand Management Program in New York~\cite{Coddington_Change_in_BQ} focuses on distribution-level capacity issues. In the planning community, DER-based approaches to long-term planning problems are often referred to as \emph{non-wire alternatives (NWAs)}. The basic premise is that NWAs can manage the load shape and peak to avoid or at least delay the need for capacity expansion.

Economically, the reason for deferring is the time-value of money, which states that a dollar spent now is (typically) more valuable than a dollar spent later~\cite{levy1994capital}. Policy-wise, the benefits of delaying capital-intensive projects are reducing the risk of the expected load not materializing and avoiding politically unpopular projects~\cite{stanton2015getting}. In this paper, we focus on the economic question and ask: \emph{is delaying traditional expansion investments worth the costs of NWAs?}

The answer to the aforementioned question is non-trivial. For one, the cost and benefits of NWAs are not only a function of their installed capacities but also of their operations. For instance, the benefits delivered by an energy storage (ES) system may include peak load reduction, load shifting to low-price hours, and reserve provision. All of these benefits depend on how the ES  system is charged and discharged, i.e., how it is operated. Thus, one must co-optimize the investment and operation of NWAs. This co-optimization can lead to a large and computationally difficult problem, especially if we consider uncertainty from renewable energy resources. Furthermore,  considering the time-value of delaying investments introduces non-linearities that result in a non-convex problem even when integer variables are not present.

Explicitly considering DERs as alternatives to capacity expansion incorporates a typically overlooked value stream to DER projects. Previous works have pointed out that delaying or avoiding capacity expansion could account for a large share of the benefits derived from DER projects (e.g., 20\% to 50\% for the case in~\cite{SIDHU2018881}). Thus, accounting for this value stream can push a DER project from ``the red'' into economic viability.

In this work, we restrict our analysis to DERs as NWAs. Non-capacity related technologies such as network management, grid optimization, data analytics, dynamic pricing, and others are outside the scope of this paper. As shown by~\cite{capitanescu2015assessing}, network reconfiguration techniques in a high distributed generation (DG)-context can serve as substitutes to network reinforcement. The results in~\cite{capitanescu2015assessing} suggests that network reconfiguration techniques, DERs, and traditional capacity expansion likely exhibit complex interactions. Considering these three elements is an interesting route for future research.

Here, we make the following contributions:
\begin{enumerate}
	\item A \emph{formulation} of the NWAs planning problem that determines i) the investment, ii) the operation of NWAs and iii) the timing of the capacity expansion. We manage load, solar generation, and energy efficiency (EE) performance uncertainty via robust optimization~\cite{houda2006comparison,bertsimas2003robust}. Our formulation explicitly includes the value of expansion delay in the objective of the problem and as an additional revenue stream for DERs. 
	\item  \emph{Tractable algorithms} for the NWAs planning problem. Modeling decades-long operation of NWAs can lead the problem to have a number of variables on the order of millions, although smaller approximations of the problem are possible (e.g., using representative days as in~\cite{CALVILLO2016340, KANDIL2018961}). Furthermore, the variables and constraints that model the timing of capacity expansion introduce non-convexities. We present two solution techniques.
	The first technique uses the \emph{Dantzig-Wolfe Decomposition Algorithm (DWDA)}. We deal with the scale of the NWAs planning problem by decomposing it into smaller subproblems. The non-convexities end up confined to a small master problem (in the order of tens to hundreds of variables) which we decompose into a small number of linear programs. In the second technique,  \emph{we fix the time of capacity expansion to eliminate the problem's non-convexities}. We sequentially solve the convex problem one time for each year in the optimization horizon and pick the best solution. 
	\item A \emph{case study} where NWAs may be used to defer substation and feeder upgrades at the University of Washington (UW) Seattle Campus using real data. We estimate the performance of the NWAs planning solution via Monte Carlo simulations.
\end{enumerate}

The idea of delaying infrastructure investment by curbing load was first introduced in~\cite{Xin_One_2005} where the authors quantify the effects of load reduction on avoided infrastructure costs. However,~\cite{Xin_One_2005} does not find the optimal amount of reduction or the appropriate technologies to do so. Similarly, the authors of~\cite{Gin_Quant_2006} and~\cite{Piccolo_Evaluating_2009} quantify the value of expansion delays by explicitly modeling DG as the mechanism to reduce net load. However, they neither address the problem of finding optimal DG investment nor consider other types of DERs. In~\cite{Samper_Investment_2013_I}, the authors determine optimal investments in DG considering the value of network investment deferral. However, their non-linear mixed-integer formulation is intractable for large systems. In contrast, we consider a broader set of NWAs and tackle the problem by solving a series of smaller convex problems.

Beyond the above-cited works, there is relatively little literature on holistic DER planning. Most consider a narrow definition of the term DER that only includes DG, e.g.,~\cite{Alarcon-Rodriguez_MO_2009, Yuan_CO_2017}, or only ES and demand response (DR)~\cite{Dvorkin_Merchant_2017}. Instead, we consider a generic definition of DERs. Furthermore, while the typical value streams of DERs include peak, energy, and greenhouse gas reduction~\cite{ GACITUA2018346, VIANA2018456}. In this paper, we explicitly consider the value of delaying capacity expansion.

Every model faces a dilemma: accuracy vs. solvability. On the one hand, an accurate model may be too complicated for the available computational resources. On the other hand, an overly simplistic model may not sufficiently reflect reality. In the context of DER planning, this trade-off is pointed out in~\cite{ALARCONRODRIGUEZ20101353}. A common technique to modulate the accuracy and solvability of a DER planning model is to modify the number of time-steps in the modeling horizon. Some change the number of representative days in a year or month (e.g., in~\cite{CALVILLO2016340, KANDIL2018961}) or modify the time step resolution (e.g., 4 hours per time step~\cite{AKBARI2016567} or a year per time step~\cite{MAHBUB20171487, WIERZBOWSKI201693}). While the relatively coarse time step resolutions in~\cite{CALVILLO2016340,KANDIL2018961, AKBARI2016567, MAHBUB20171487, WIERZBOWSKI201693} may be appropriate for their specific studies, there are two issues with the approach. First, the practitioner needs to know ``how coarse is too coarse'' and/or how to select representative days. Second, one may want to consider DER technologies whose modeling require high time-step resolutions\footnote{For instance, relatively high time-step resolutions are needed to adequately capture inter-temporal dynamics relevant to ES operation, as tacitly noted in~\cite{Deane_2013_Derivation}. DR is another technology that requires relative high time-step resolution for realistic, efficient, and reliable modeling and operation~\cite{Neves_2015_demand}.  }. 

We avoid grappling with the accuracy vs. solvability dilemma by adopting a relatively high-resolution time step (1 hour over a 20-year horizon in our case study) and achieving solvability by decomposing the problem. We take this approach for a couple of reasons. The first is that we attempt to develop a method that applies to a wide range of current and future DER technologies. As mentioned previously using high time resolutions is crucial for technologies like ES and DR. The second reason for modeling the planning horizon using high-resolution time steps is that we avoid the expertise needed to implement reduction and simplification techniques. Thus, we avoid the need of additional algorithms that select ``representative’’ days, months, or weeks as done in works like~\cite{CALVILLO2016340, KANDIL2018961}.

In addition to being sufficiently accurate and solvable, implementation of a model must be practical. That is, the economic, time, and human costs of implementation and data collection must be reasonable. In our case, practicality means that our method demands a reasonable amount of engineering knowledge, data quality, computing power, and other relevant resources from the practitioner (e.g., a utility). While the method presented in this paper is mathematically complex by traditional utility's standards, in the mathematical model would run in the background in commercial deployments of our method. The end user would not necessarily need to make ``under-the-hood'' changes to the math behind our model. We believe that most utilities in North America can utilize our method with reasonable to little additions to the engineering skills and computational resources that they already possess in their planning departments. In fact, some leading utilities are already experienced in NWAs planning and posses a skilled workforce and resources to easily implement our method. Regarding input data, our method does not require any data in addition to the one required in more conventional DER planning methods.

This paper differs and adds to our previous work published in~\cite{Contreras-Ocana_2018_non} in the following ways. 
\begin{itemize}
	\item We now consider uncertainty and present a robust formulation to manage it. 
	\item Furthermore, we adopt the scenario generation approach presented in~\cite{chen2018model} to produce the load profiles used to formulate the robust NWAs planning problem. 
	\item We provide a probabilistic assessment of the solution to the NWAs using Monte Carlo simulations.
	\item We improve the solution techniques and provide results on computational performance.
\end{itemize}

This paper is organized as follows. Section~\ref{sec:Substation_upgrade_problem} states the capacity expansion problem and formulates it as an optimization problem.  Section~\ref{sec:NWA} formulates a generic NWA model and four specific technologies: EE, PV, DR, and ES. Section~\ref{sec:NWA_PP} describes the NWAs planning problem and uncertainty modeling. Section~\ref{sec:DWDA} provides two solution techniques for the planning problem. Section~\ref{sec:case_study} presents a case study of load-growth and NWAs at the UW. Section~\ref{sec:conclusion} concludes the paper and provides suggestion for future research.

\section{The capacity expansion problem}
\label{sec:Substation_upgrade_problem}

System planners typically like to expand capacity at the \emph{latest}\footnote{In practice, planners act with a lead or cushion time to account for implementation uncertainties and delays.} possible time that meets the expected load growth~\cite{SCL_planning_practices_2018}. The main economic reason to delay capacity expansion is the time-value of money: we would like to spend a dollar later rather than now. Let $l_a^\mathrm{p}$ denote expected peak load\footnote{The expected peak load $\boldsymbol{l}^\mathrm{p}$ is typically forecasted using a variety of inputs such as population and economic growth projections, planned construction projects, weather forecasts, etc~\cite{seifi2011electric}.} during year $a$ and the pre-expansion capacity as $\overline{l}$. The vector of expected peaks in the planning horizon is denoted as $\boldsymbol{l}^\mathrm{p}$.  After expansion, we assume that any reasonable load can be accommodated during the planning horizon.  Then, the decision rule for choosing a year to expand capacity is
\begin{equation} \label{eq:subs_upgrade_prob}
\mathrm{CapEx}(\boldsymbol{l}^\mathrm{p})=\delta\; | \;l^\mathrm{p}_{\delta+1}>\overline{l}, \;l_a^\mathrm{p} \le \overline{l}  \;\forall \; a \le \delta.
\end{equation}
The decision rule $\mathrm{CapEx}$ states that the planner expands capacity at a future year $\delta$, immediately before the limit $\overline{l}$ is \emph{first} reached by the load. In practice, an expansion project for year $a$ should be be started with enough cushion time to account for implementation uncertainty and delays. 

In practice, systems have different types of capacity limits such as transient, steady-state, and/or security-induced limits (e.g., the N-1 security criterion). In this paper, we ignore transients and focus on the steady-state capacity limit that meets security requirements. Additionally, we sidestep non-capacity technical questions such as voltage stability~\cite{kundur1994power} and neglect the losses in the system.

Furthermore, we analyze capacity expansion at a single point in a distribution system (e.g., a feeder or substation) and assume that the downstream network is not congested. Consequently, the problem of determining the network placement of each NWA is outside the scope of this work. Nevertheless, we acknowledge the importance of the NWAs placement problem as it has an important influence on the economic viability of a DER project~\cite{casten2002optimizing}. We believe that these are reasonable assumptions because:
\begin{enumerate}
	\item We envision the proposed framework to be among the initial stages of a NWAs project and is therefore mainly focused on capacity. Phenomena such as voltage stability are of definite importance in system operations, but these can be addressed after an expansion plan is decided (e.g., by using DERs to provide voltage support~\cite{kashem2005multiple}, voltage regulators and/or, capacitor banks~\cite{kundur1994power}). However, a capacity shortage is harder to be addressed afterward and must be considered in the early planning stages. Losses impact the economics of the system, but is secondary during the capacity planning stage~\cite{wang1994modern}.
	\item  Since most distribution systems are (approximately) radial~\cite{Kersting01}, if a point in the downstream network is congested, the framework proposed in this paper can be applied to the subnetwork. Furthermore, because NWAs tend to reduce net peak demand (e.g., PV and EE offset load while DR and ES by shifting peak load), congestion is less likely than in a scenario where NWAs are not present.
\end{enumerate}

Let $I$ denote the inflation-adjusted cost of capacity expansion. Here, we assume that the inflation-adjusted cost of capacity expansion is constant throughout the planning horizon. Then, if system capacity is expanded at year $\mathrm{CapEx}(\boldsymbol{l}^\mathrm{p})$, the present cost of the investment is
\begin{equation}
\tilde{I}(\boldsymbol{l}^\mathrm{p}) = \frac{I} {\left(1+\rho\right)^{\mathrm{CapEx}(\boldsymbol{l}^\mathrm{p})}}  \label{eq:IW_PV}
\end{equation}
where $\rho$ is the annual discount rate, a positive quantity closely related to the interest rate~\cite{seifi2011electric,levy1994capital}.

Our goal would be to minimize $\tilde{I}(\boldsymbol{l}^\mathrm{p})$ via the variable $\mathrm{CapEx}(\boldsymbol{l}^\mathrm{p})$. However, since $\mathrm{CapEx}(\boldsymbol{l}^\mathrm{p})$ itself is given by a decision rule in~\eqref{eq:subs_upgrade_prob}, it is more convenient to write the present cost $\tilde{I}(\bd{l}^\mathrm{p})$ as an minimization problem. Let the planning horizon be denoted as $\bd {\mathcal A}$. For example, $\bd{\mathcal{A}}=\{1,2,\dots,20\}$ is a 20 year planning horizon. We have the following lemma.
\begin{lemma} \label{thm:IW_PV_CVX}
	The function $\tilde{I}$ from~\eqref{eq:IW_PV} can be reformulated as the optimization problem
	\begin{subequations}  \label{eq:IW_PV_CVX}
		\begin{align}
		\tilde{I}(\bd{l}^\mathrm{p})=\min_{\delta}& \frac{I}{(1+\rho)^{\delta}}  \\
		\mbox{s.t. } & 0 \le \delta \le |\bd {\mathcal A}| \\
		& l_{a}^\mathrm{p} \le \overline{l}\;\; \forall \;a<\delta. \label{eq:NCVX_Const}
		\end{align}
	\end{subequations}
\end{lemma}
The reformulation of $\tilde{I}$ allows us to embedded it in an optimization problem without the need of including conditionals in~\eqref{eq:subs_upgrade_prob}. The proof of Lemma~\ref{thm:IW_PV_CVX} is given in Appendix~\ref{sec:app_proof_lemma}. 

Note that for any given set of yearly peak loads $\bd{l}^\mathrm{p}$,~\eqref{eq:IW_PV_CVX} is convex. However, if we treat the peak load as a function of NWAs operation, and therefore as an optimization variable,~\eqref{eq:IW_PV_CVX} becomes non-convex.

While it is unfortunate that~\eqref{eq:IW_PV_CVX} is non-convex in $\bd{l}^\mathrm{p}$, Section~\ref{sec:DWDA} shows how to handle its non-convexities by solving at most $|\bd{ \mathcal A}|$ small-scale linear problems. This is because the length of the planning horizon, $|\bd {\mathcal A}|$, is in the order tens of years. These small problems can be reliably solved using off-the-shelf solvers (e.g., Gurobi). In the next section, we introduce the generic model of an NWA and instantiate it.

\section{Non-wire alternatives}
\label{sec:NWA}

\begin{table}[]
	\centering
	\caption{Nomenclature for Section~\ref{sec:NWA}}
	\label{table:Nomenclature}
	\begin{tabular}{@{}|l|l|@{}}
		\toprule
		\rowcolor[HTML]{EFEFEF} 
		\multicolumn{2}{|l|}{\cellcolor[HTML]{EFEFEF}\textbf{Sets}}                          \\ \midrule
		$\bd{\mathcal{A}}$         & Set of years in the planning horizon, indexed by $a$                         \\ \midrule
		$B^\mathrm{EE}$         & Number of EE segments, indexed by $b$                           \\ \midrule
		$\bd{\mathcal{N}}$         & Set of NWA technologies, indexed by $i$                           \\ \midrule
		$\bd{\mathcal{T}}$         & Set of operating intervals in one year, indexed by $t$                       \\ \midrule
		$\bd{\mathcal{X}}_i$         & Set of feasible operating regimes of NWA $i$                        \\ \midrule
		$\bd \Phi_i$          &Set of feasible investment decisions of NWA $i$                        \\ \midrule
		\rowcolor[HTML]{EFEFEF} 
		\multicolumn{2}{|l|}{\cellcolor[HTML]{EFEFEF}\textbf{Decision variables}}                          \\ \midrule
		$c_{a,t}$/$d_{a,t}$           & ES charge/discharge at $a,t$ (MW)      \\ \midrule
		$g_{a,t}^\mathrm{PV}$          & PV generation at $a,t$  (MW)                          \\ \midrule
		$g^\mathrm{PV}_\mathrm{CAP}$          & PV installed capacity (MW)                           \\ \midrule
		$I^\mathrm{NW}_i$          & Investment cost of NWA $i$  (\$)               \\ \midrule
		$\tilde{I}^\mathrm{NW}_i$          & Present cost of capacity expansion  (\$)               \\ \midrule
		$l^i_{a,t}$ / $l_{a,t}$          &  Load of NWA $i$ / total load at $a,t$ (MW)                   \\ \midrule
		$l^\mathrm{p}_{a}$          & Peak load during year $a$ (MW)                   \\ \midrule
		$r^\mathrm{a,t}_\mathrm{CAP}$          &  Load reduction at $a,t$ from DR deployment (MW)                           \\ \midrule
		$r^\mathrm{DR}_\mathrm{CAP}$          & Capacity of DR-enabled load (MW)                           \\ \midrule
		$r^\mathrm{EE}_{a,t}$          & Load reduction from EE at $a,t$  (MW)                           \\ \midrule            
		$s_{a,t}$           & ES state-of-charge at $a,t$ (MWh)      \\ \midrule
		$s^\mathrm{max}_a$          & \shortstack[l]{ES capacity at year $a$  (MWh) \\ Note: $s^\mathrm{max}_0$ is the initial capacity}                         \\ \midrule
		$\bd x_i$          & Operating decision variables of NWA $i$                                \\ \midrule
		$\epsilon_b^\mathrm{EE}$    & Percentage load reduction of EE segment $b$ (\%)                             \\ \midrule
		$\bd \phi_i$          & Investment decision variables of NWA $i$                                 \\ \midrule
		\rowcolor[HTML]{EFEFEF} 
		\multicolumn{2}{|l|}{\cellcolor[HTML]{EFEFEF}\textbf{Parameters}}                           \\ \midrule
		$C_b^\mathrm{EE}$                     & Cost of segment $b$ EE $\mathrm{\left(\frac{\$}{\%}\right)}$              \\ \midrule
		$C^\mathrm{ES}$/$C^\mathrm{DR}$                      & ES/DR cost $\mathrm{\left(\frac{\$}{MW}\right)}$            \\ \midrule        
		$C^\mathrm{PV}$                     & PV cost $\mathrm{\left(\frac{\$}{MW}\right)}$            \\ \midrule    
		$\overline{g}^\mathrm{PV}_\mathrm{CAP}$          & Maximum PV installed capacity (MW)                           \\ \midrule        
		$l_{a,t}^\mathrm{b}$                      & Base load at $a,t$  (MW)              \\ \midrule
		$\overline{s}^\mathrm{max}_0$      & Maximum ES capacity $t$  (MWh)              \\ \midrule
		$\alpha^\mathrm{DR}$                      & DR rebound coefficient                      \\ \midrule
		$\alpha^\mathrm{EPR}$                      & Energy-to-power ratio of the ES system                      \\ \midrule
		$\alpha_{a,t}^\mathrm{EE}$                      & Accuracy of projected load reduction at $a,t$ (\%)                     \\ \midrule
		$\alpha_{a,t}^\mathrm{PV}$                      & Per-unit PV generation at $a,t$ (\%)                   \\ \midrule
		$\beta^\mathrm{ESD}$    & ES degradation coefficient $\mathrm{\left(\frac{MWh}{MW}\right)}$             \\ \midrule
		$\Delta t$                      & Length of operating time interval (hours)                        \\ \midrule
		$\overline{\epsilon}_b^\mathrm{EE}$                      & Size of EE segment $b$ (\%)                       \\ \midrule
		$\eta_c$/$\eta_d$ & ES charge/discharge efficiency                        \\  
		\bottomrule
	\end{tabular}
\end{table}

Let the index $i$ denote a NWA technology. A generic NWA is characterized by six elements:
\begin{enumerate}[noitemsep,topsep=0pt]
	\item investment (or sizing) decision variables $\bd \phi_i$, for example, the energy capacity of an ES system;
	\item operating decision variables $\bd x_i$, for example, ES hourly charging and discharging decisions;
	
	\item a set of feasible investment decisions $\bd \Phi_i$, for example, the set of ES systems that physically fit in a site;
	
	\item a set of feasible operating regimes $\bd{\mathcal{X}}_i(\bd \phi_i)$, for example, the set of hourly charging and discharging decisions that comply with charge, discharge, and state-of-charge limits;
	
	\item a set of functions $l^i_{a,t}(\bd x_i)$ that map operating decisions onto load at time $t$ of year $a$, for example, for an ES system at time $a,t$ the load is defined as charge minus discharge; and
	
	\item an investment cost function $I^\mathrm{NW}_i(\bd \phi_i)$, for example, the investment cost of an ES system is the energy capacity times the per-kWh cost of storage.
	
\end{enumerate}

While investment decisions are made once (or at most a few times) in a planning horizon, operating decisions are made frequently. In this paper, the operating decision time intervals length is $\Delta t$ hours and $\bd{\mathcal{T}}$ denotes the set of operating intervals in one year.

It is worth noting that our framework allows for the consideration of NWAs that are installed in the future. This is especially important since DERs are expected to have significant technology and cost-related developments. For instance, it is possible to consider installing an ES system in during year $a$ with future technology and cost characteristics. To do so, one would consider an ES technology that starts operating at year $a$ and whose cost is time-value adjusted.

We assume that $\bd \Phi_i$ and $\bd{\mathcal{X}}_i(\bd \phi_i)$ are convex, $I^\mathrm{NW}_i(\phi_i)$ is a convex function, and that the functions $l^i_{a,t}(\bd x_i)$ are linear in $\bd x_i$. While these assumptions restrict the complex reality of NWAs to simpler models, they provide computational tractability in an optimization context. The restrictions and computational benefits implied by the aforementioned assumptions are accepted and widely utilized in the academic literature and by practitioners, e.g., in~\cite{SCL_planning_practices_2018, Shi_2018_Using, Xu_2017_Scalable, Gantz_2014_Optimal, Nasir_2018_Optimal, Muneer_2011_Large, Alharbi_2018_Stochastic, Carpinelli_2013_Optimal, Dvorkin_2018_Co-Planning}. Now, we describe each of the six elements that characterize a NWA for the four technologies that we consider in this paper: EE, PV, DR, and ES.

\paragraph*{Energy efficiency~(EE)} For EE, the investment decision is to choose a percentage base load reduction that translates into a $r_{a,t}^\mathrm{EE}$ reduction at every time period. We model the investment cost, $I^\mathrm{NW}_\mathrm{EE}$, as a convex piece-wise linear function of load reduction~\cite{brown2008us}. The slope of each of the $B^\mathrm{EE}$ segments, $C_b^\mathrm{EE}$, represents the marginal cost of load reduction. The six parameters that define EE as a NWA are
\begin{align*}
	&\bd \phi_\mathrm{EE} = \left\{ \epsilon_b^\mathrm{EE}\right\}_{b=1,\hdots,B^\mathrm{EE}} ,\;\; \bd x_\mathrm{EE} = \{r^\mathrm{EE}_{a,t}\}_{a\in \bd {\mathcal A},\; t \in \bd{\mathcal{T}}},\\
	&\bd \Phi_\mathrm{EE} \! =\left\{ \epsilon_b^\mathrm{EE} \; | \;  \epsilon_b^\mathrm{EE}\in \left[ 0, \overline{\epsilon}_b^\mathrm{EE}\right] \; \forall\; b=1\hdots, B^\mathrm{EE} \right\} ,\\
	&\bd{\mathcal{X}}_\mathrm{EE} \!\left(\bd \phi_\mathrm{EE}\right) \! = \! \left\{r^\mathrm{EE}_{a,t}  | r^\mathrm{EE}_{a,t} = \alpha_{a,t}^\mathrm{EE} \cdot  l_{0,t}^\mathrm{b} \!\cdot\!\! \sum_{b=1}^{B^\mathrm{EE}}\! \epsilon^{EE}_b \forall  a \in  \bd {\mathcal A},t \in \bd{\mathcal{T}}\right \},\\
	&l_{a,t}^\mathrm{EE}(\bd x_\mathrm{EE}) = -r_{a,t}^\mathrm{EE} ,\;\;I^\mathrm{NW}_\mathrm{EE}(\bd \phi_\mathrm{EE}) =   \sum_{b=1}^\mathrm{B^\mathrm{EE}} C^\mathrm{EE}_b \cdot \epsilon^\mathrm{EE}_b ,
\end{align*}
where $\epsilon_b^\mathrm{EE}$ is the projected percentage reduction for each piece-wise linear segment of $I^\mathrm{NW}_\mathrm{EE}$, $\overline{\epsilon}_b^\mathrm{EE}$ is the size of each segment, and $l_{0,t}^\mathrm{base}$ is the base load (i.e., the pre-EE load). The parameter $\alpha_{a,t}^\mathrm{EE}$ represents accuracy the projected load reduction. For instance, $\alpha_{a,t}^\mathrm{EE}=1$ represents no error while $\alpha_{a,t}^\mathrm{EE}=0.9$ represents a 10\% underestimation.

Note that more complex models of EE are allowed within our framework. For instance, each EE reduction segment $b$ could be associated with a load type (e.g., a segment dedicated to heating, ventilation, and air conditioning (HVAC) and a segment dedicated to artificial lighting) and specific costs (e.g., HVAC retrofit costs and costs of upgrading artificial lighting to LEDs). Thus, implementing an HVAC retrofit would translate into a load reduction during time periods where HVAC is in use. Similarly, upgrading artificial lighting to LEDs would translate into a load reduction during time periods when artificial lighting is in use.

\paragraph*{Solar photovoltaic generation~(PV)}
The PV investment decision is the installed capacity $g^\mathrm{PV}_\mathrm{CAP}$.  The solar energy generation at time $t$, $a$ is $g^\mathrm{PV}_{a,t} = \alpha_{a,t}^\mathrm{PV}\cdot g^\mathrm{PV}_\mathrm{CAP}$ where $\alpha_{a,t}^\mathrm{PV} \in [0,1]$ is solar generation per unit of PV installed capacity and is related to solar irradiation. For instance, at night, when solar irradiation is zero, $\alpha_{a,t}^\mathrm{PV}=0$ and when PV system outputs its full capacity, $\alpha_{a,t}^\mathrm{PV}=1$. The parameters that define solar PV as a NWA are
\begin{align*}
	& \bd \phi_\mathrm{PV}  =  g^\mathrm{PV}_\mathrm{CAP} ,\;\; \bd x_\mathrm{PV} = \{g_{a,t}^\mathrm{PV}\}_{a\in \bd {\mathcal A},\; t \in \bd{\mathcal{T}}}, \\
	&\bd \Phi_\mathrm{PV} = \left\{ g^\mathrm{PV}_\mathrm{CAP}\; |\; g^\mathrm{PV}_\mathrm{CAP} \in \left[0, \overline{g}^\mathrm{PV}_\mathrm{CAP} \right] \right\}, \\
	& \bd{\mathcal{X}}_\mathrm{PV} \left(\bd \phi_\mathrm{PV}\right) = \{g_{a,t}^\mathrm{PV} \; |\;g_{a,t}^\mathrm{PV} = \alpha_{a,t}^\mathrm{PV}\cdot  g^\mathrm{PV}_\mathrm{CAP}  \;\forall \; a\in \bd {\mathcal A},\; t \in \bd{\mathcal{T}}\} ,\\
	& l_{a,t}^\mathrm{PV}(\bd x_\mathrm{PV}) = -g_{a,t}^\mathrm{PV},\;\;I^\mathrm{NW}_\mathrm{PV}(\bd \phi_\mathrm{PV}) = C^\mathrm{PV} \cdot g^\mathrm{PV}_\mathrm{CAP},
\end{align*}
where $C^\mathrm{PV}$ is the per-unit PV capacity cost  and $\overline{g}^\mathrm{PV}_\mathrm{CAP}$ is the PV capacity limit.

\paragraph*{Demand response~(DR)}
We consider investments in DR communication and control infrastructure that enable shifting a portion of the load. The investment decision is the amount DR-enabled load $r^\mathrm{DR}_\mathrm{CAP}$ which limits the demand reduction $r_{a,t}^\mathrm{DR}$ deployed during year $a$, operating period $t$. A load reduction $r_{a,t}^\mathrm{DR}$ causes a demand rebound of $\alpha^\mathrm{DR}\cdot r_{a,t}^\mathrm{DR}$ during time period  $t+1$. The coefficient $\alpha^\mathrm{DR}$ is a number $\ge 1$ and is related to efficiency losses caused by DR deployment~\cite{lutolfimpact}. More sophisticated rebound models such as the ones in~\cite{lutolfimpact} are admissible in our framework. The parameters that define DR are
\begin{align*}
	&\bd \phi_\mathrm{DR}  =  r^\mathrm{DR}_\mathrm{CAP} ,\;\; \bd x_\mathrm{DR} = \{r_{a,t}^\mathrm{DR}\}_{a\in \bd {\mathcal A},\; t \in \bd{\mathcal{T}}},\\
	& \bd \Phi_\mathrm{DR} = \left\{ r^\mathrm{DR}_\mathrm{CAP}\; | \; r^\mathrm{DR}_\mathrm{CAP} \in\left[ 0, \overline{r}^\mathrm{DR}_\mathrm{CAP}  \right] \right\},\\
	&\bd{\mathcal{X}}_\mathrm{DR} \left(\bd \phi_\mathrm{DR}\right) = \{ r_{a,t}^\mathrm{DR} |  r_{a,t}^\mathrm{DR} \in\left[ 0, r^\mathrm{DR}_\mathrm{CAP} \right] \; \forall \; a\in \bd {\mathcal A},\;t \in \bd{\mathcal{T}} \} ,\\
	& l_{a,t}^\mathrm{DR}(\bd x_\mathrm{DR}) =\alpha^\mathrm{DR}  \cdot r_{a,t-1}^\mathrm{DR} -r_{a,t}^\mathrm{DR} ,\;\; I^\mathrm{NW}_\mathrm{DR}(\bd \phi_\mathrm{DR}) = C^\mathrm{DR}  \cdot r^\mathrm{DR}_\mathrm{CAP},
\end{align*}
where $C^\mathrm{DR}$ is the per-unit cost of DR. In our work, we ignore binary variables that arise from fixed DR or customer enrollment costs.

In practice, DR-enabled load is limited by a number of technical, regulatory, and human factors. We summarize these limitations by limiting DR installed capacity to $\overline{r}^\textrm{DR}_\textrm{CAP}$. It is important to note that the model could be refined if the planner possesses  better knowledge of the DR capabilities of the load.

\paragraph*{Lithium-ion energy storage}
The ES investment decision is the initial (e.g., name-plate capacity) energy capacity $s^\mathrm{max}_0$ of the storage system. The operating variables are the charge $c_{a,t}$, discharge $d_{a,t}$ and the state-of-charge $s_{a,t}$. In addition, we consider  energy capacity as an operating variable since $s^\mathrm{max}_a$ may be different than $s^\mathrm{max}_0$ because we model battery degradation. The feasible operating region of the ES system is defined by~\eqref{eq:feas_reg_start}-\eqref{eq:feas_reg_end} and include the usual charge, discharge, and state-of-charge limits~\cite{Sarker_Optimal_2017}. An important aspect of chemistry-based ES that should be taken into account in planning algorithms is that its lifetime is related to its operation.  As expressed in Eq.~\eqref{eq:degradation}, the storage capacity degrades by $\beta^\mathrm{ESD}$ per-unit charge/discharge~\cite{Sarker_Optimal_2017}. Our approach to degradation modeling is in line with the methods used in~\cite{Ru_2013_Storage, Atia_2016_Sizing}. The parameters that define ES are
\begin{subequations}
	\begin{align}
		&\bd \phi_\mathrm{ES}  =   s^\mathrm{max}_0 ,\; \bd x_\mathrm{ES} = \{c_{a,t},\; d_{a,t},\; s_{a,t},\; s_a^\mathrm{max}\}_{a\in \bd{\mathcal A},\; t \in \bd{\mathcal{T}}},\\
		&\bd \Phi_\mathrm{ES} =  \left\{  s^\mathrm{max}_0 \;|\; s^\mathrm{max}_0 \in \left[0,\overline{s}^\mathrm{max}_0 \right] \right\},\\
		&\bd{\mathcal{X}}_\mathrm{ES} \left(\bd \phi_\mathrm{ES}\right) = \biggl\{c_{a,t},\; d_{a,t},\; s_{a,t},\; s_a^\mathrm{max}\;| \label{eq:feas_reg_start} \\
		& s_{a,t+1} = s_{a,t} +  \Delta t\cdot \left( \eta_c \cdot  c_{a,t} -\frac{d_{a,t}}{\eta_d} \right) \; \forall \; a\in\bd {\mathcal A},\; t \in \bd{\mathcal{T}}, \\
		&s_{a}^\mathrm{max} = s^\mathrm{max}_0  - \beta^\mathrm{ESD} \cdot\sum_{k = 1}^{a-1} \sum_{t \in \bd{\mathcal{T}}} (c_{k,t} + d_{k,t}) \; \forall \; a\in\bd {\mathcal A} , \label{eq:degradation} \\
		& s_{a,t} \in \left[0, s_a^\mathrm{max}\right],\;c_{a,t} ,d_{a,t} \in \left[0, \frac{s^\mathrm{max}}{\alpha^\mathrm{EPR}}\right] \; \forall \; a\in \bd {\mathcal A},\; t \in \bd{\mathcal{T}}  \biggr\}, \label{eq:feas_reg_end} \\
		&l_{a,t}^\mathrm{ES}(\bd x_\mathrm{ES}) =c_{a,t} -d_{a,t} ,\;\; I^\mathrm{NW}_\mathrm{ES}(\bd \phi_\mathrm{ES}) = C^\mathrm{ES} \cdot s^\mathrm{max}_0,
	\end{align}
\end{subequations}
where $\eta_c \; (\eta_d)$ is the charge (discharge) efficiency, $\alpha^\mathrm{EPR}$ is the energy-to-power ratio of the ES system, and  $C^\mathrm{ES}$ is the dollar per-unit energy cost of ES capacity.  We consider investments in lithium-ion ES because of their ubiquity although other chemistries are compatible with our framework.

In our work, we model the investment cost of each NWA as a linear function of installed capacity. In practice, investment costs functions depend on a multiplicity of factors and are not necessarily linear. Some components of the investment cost that are not necessarily linear to capacity depend on labor, permitting, overhead, and inverter costs~\cite{Fu_2018_Evaluating}. However, we believe that adopting a linear investment cost function is an appropriate simplification that help keeps the problem computationally manageable. Similar investment cost functions are widely adopted in the literature, e.g., in~\cite{Nasir_2018_Optimal, Muneer_2011_Large, Alharbi_2018_Stochastic, Carpinelli_2013_Optimal}.

Our framework allows for other NWAs to be included, for example, electric vehicles, a diverse range of ES technologies, dispatchable DG, etc. In this paper, we limit our consideration to the four technologies described above (EE, PV, DR, ES) because they are the most mature, do not produce carbon emissions, and readily deployable in an urban environment, which fits our case study about the UW.

\subsection*{The capacity expansion problem in a NWAs context}

Let the total load including a set $\bd{\mathcal{N}}$ of NWAs be denoted by $l_{a,t}(\bd{x}) = l^\mathrm{b}_{a,t} + \sum_{i \in \bd{\mathcal{N}}} l_{a,t}^i(\bd{x}_i)$ where $\bd{x} = \{\bd x_i \}_{i\in \bd{\mathcal{N}}}$. Then, the yearly peak load as a function of NWA operation is $l_a^\mathrm{p}(\bd{x}) = \max_{t \in \bd{\mathcal{T}}}\{ l_{a,t}(\bd{x}) \}$ where $\bd x$ denote the operating decisions of all the NWAs.

The only decision in the traditional capacity expansion problem is the time of capacity expansion. With NWAs, however, the present cost of expansion
\begin{equation}
	\tilde{I}(\bd{l}^\mathrm{p}(\bd{x}))  \label{eq:NWAs_planning_PC_expansion}
\end{equation}
is a function of the NWAs operation and gives the planner the opportunity to invest in and operate a set of NWAs that minimizes~\eqref{eq:NWAs_planning_PC_expansion}. However, a good plan should consider additional NWAs costs and benefits (e.g., demand charge reductions, DR rebound costs, etc.). In the next section, we present a holistic NWAs planning problem that decides the investment and operation of NWAs, and the time of capacity expansion.

\section{The non-wire alternatives planning problem}
\label{sec:NWA_PP}
\subsection{Deterministic formulation}

	The NWAs planning problem
	\begin{equation}
	\min_{\substack{\bd \phi_i \in \bd \Phi_i \\ \bd x_i \in \bd{\mathcal{X}}_i (\bd \phi_i) }} \!\!\biggl\{\sum_{i \in \bd{\mathcal{N}}} \!\! \left[ C_i^\mathrm{O}(\bd x_i)\! + \!I_i^\mathrm{NW}(\bd \phi_i) \right]  \!+C^\mathrm{D}(\boldsymbol{l}^\mathrm{p}(\boldsymbol{x}))  + \tilde{I} (\boldsymbol{l}^\mathrm{p}(\boldsymbol{x})) \biggr\} \label{prob:NWA_planning}
	\end{equation}
	minimizes NWA operating costs (e.g., fuel costs), NWA investment costs, peak load-related costs (e.g., peak-demand charge), the present cost of capacity expansion and maximizes operating benefits (e.g., ES shifting load to low-cost hours or remuneration for services to the market). The NWA operating costs and benefits are expressed by $C_i^\mathrm{O}(\bd x_i)$ and include, if applicable, secondary value streams such as revenue from providing reserves, frequency regulation, or backup capacity. The investment costs for NWA $i$ are denoted by $I_i^\mathrm{NW}(\bd \phi_i)$ and peak load-related costs by $C^\mathrm{D}(\boldsymbol{l}^\mathrm{p}(\boldsymbol{x}))$. The objective of~\eqref{prob:NWA_planning} is subject to investment and operating constraints of each NWA.
	
	Notice that the formulation of~\eqref{prob:NWA_planning} co-optimizes a number of DER costs, benefits, and secondary value streams. As suggested by Shi et al. in~\cite{Shi_2018_Using} co-optimizing services can deliver superlinear benefits, i.e., the total co-optimized benefits are larger than the sum of individually-optimized benefits. The co-optimization approach stands in contrast to partitioning approaches such as the one proposed by Gantz et al. in~\cite{Gantz_2014_Optimal} where a portion of the DER (ES in their particular case) capacity is reserved for each service.

It is worth noting that we assume that a single stakeholders (e.g., the customer) is responsible for all the costs in~\eqref{prob:NWA_planning}. Or equivalently, we assume that, in the case of multiple stakeholders, all parties are interested in the lowest-cost solution\footnote{In practice, there might be conflicting interests among stakeholders. We believe that even in such cases, the least-cost solution is of interest (e.g., for benchmarking). Treating conflicting interests is beyond the scope of our work.}.

From~\eqref{eq:IW_PV}, $\tilde I$ contains the condition-based function $\mathrm{CapEx}$ which makes incorporating $\tilde I$ in large-scale optimization problems difficult and Problem~\eqref{prob:NWA_planning} intractable. Using Lemma~\ref{thm:IW_PV_CVX}, Theorem~\ref{thm:NWAPP_cvx} shows a more convenient formulation of the planning problem.
	\begin{thm} \label{thm:NWAPP_cvx}
		Problem~\eqref{prob:NWA_planning} is equivalent to
		\begin{subequations}
			\begin{align}
			\min_{\substack{ \bd x_i, \;\bd \phi_i ,\; \delta \\ \boldsymbol{l}^\mathrm{p}=\{l_a^\mathrm{p}\}_{a\in\bd {\mathcal A}} }} &   \sum_{i \in \bd{\mathcal{N}}} \!\left[C_i^\mathrm{O}(\bd x_i) + I_i^\mathrm{NW}(\bd \phi_i) \right]
			+C^\mathrm{D}(\boldsymbol{l}^\mathrm{p})+ \frac{I}{(1 +\rho)^{\delta}}  \\
			\mathrm{s.t. } & \bd \phi_i \in \bd \Phi_i \; \forall \; i \in \bd{\mathcal{N}}\\
			&\bd x_i \in \bd{\mathcal{X}}_i (\bd \phi_i) \; \forall \; i \in \bd{\mathcal{N}} \\
			& l^\mathrm{b}_{a,t} + \sum_{i \in \bd{\mathcal{N}}} l_{a,t}^i(\bd x_i)  \le l_a^\mathrm{p} \;\forall \; a \in \bd {\mathcal A},\; t\in\bd{\mathcal{T}} \label{eq:CVX_problem_c3}\\
			& l_{a}^\mathrm{p} \le \overline{l}\; \forall \; a<\delta \label{eq:CVX_problem_c4} \\
			& 0 \le \delta \le |\bd {\mathcal A}|. \label{eq:CVX_problem_c4.1}
			\end{align} \label{eq:CVX_problem}
		\end{subequations}
	\end{thm}
The proof of Theorem~\ref{thm:NWAPP_cvx} is given in Appendix~\ref{sec:appendix_proof_thm}. The objective of~\eqref{eq:CVX_problem} is convex because we assume that $C_i^\mathrm{O}(x_i)$, $I_i^\mathrm{NW}(\phi_i)$ and $C^\mathrm{D}(\boldsymbol{l}^\mathrm{p})$ are convex, and $\frac{I}{(1 +\rho)^{\delta}}$ is also convex. Constraint~\eqref{eq:CVX_problem_c4}, however, introduces non-convexities to the feasible solution space. Section~\ref{sec:DWDA} shows how we decompose~\eqref{eq:CVX_problem} and deal with the large-scale and non-convex nature of the problem. In the rest of this section, we show how to treat uncertainties in the planning problem.

\subsection{Uncertainty modeling} \label{subsec:uncertainty model}
We consider three major sources of uncertainty: solar irradiation ($\alpha_{a,t}^\mathrm{PV}$), base load ($l_{a,t}^\mathrm{b}$), and projected load reduction from EE measures ($\alpha_{a,t}^\mathrm{EE}$). We formulate the NWAs planning problem as a robust problem for three main reasons. First, it is less computationally-intensive than its stochastic counterpart. Second, it does not require forecasting the density functions of the uncertain parameters\footnote{We only need to estimate the maximum and minimum possible values.}~\cite{houda2006comparison,bertsimas2003robust}. Lastly and perhaps most importantly, utility planning practices typically focus on worst-case realizations. Thus, a robust approach to NWAs planning is likely more attractive to electric utilities. We direct the interested reader to Appendix~\ref{app:robust} for details on the robust formulation.

Although our robust formulation does not need scenarios of the uncertain parameters, we use them for two reasons. First, the scenarios allow us to estimate the maximum and minimum values of the uncertain parameters that are required to formulate the robust problem. And second, it allows us to evaluate the performance of the NWAs planning solution (e.g., via Monte Carlo simulation as in the case study in this paper).

	\begin{figure} 
		\centering
		\includegraphics[width=0.5\textwidth]{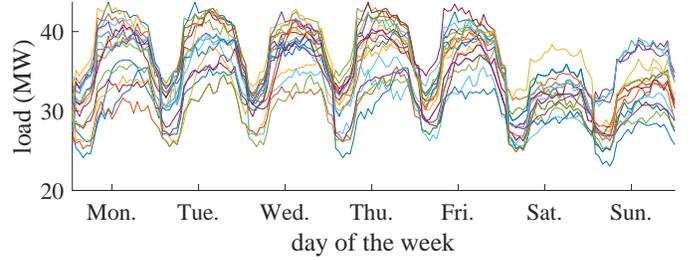}
		\caption{Twenty weekly load profiles produced via the GANs-based scenario generation technique proposed in~\cite{chen2018model}.   }\label{fig:load_scenarios_week}
	\end{figure}
	
	We adopt the scenario generation technique introduced in~\cite{chen2018model} to produce load profile and solar irradiation scenarios. The technique relies on Generative Adversarial Networks~(GANs), a machine learning-based generative model~\cite{goodfellow2014generative}. We base the GANs on a game theory setup between two deep neural networks, the \emph{generator} and the \emph{discriminator}. The generator transforms input from a known distribution $\mathbb{P}_Z$~(e.g., Gaussian) to an output distribution $\mathbb{P}_{G}$. On the other hand, the discriminator discerns historical data $\mathbb{P}_{X}$ from the output distribution $\mathbb{P}_{G}$. 
	
	In the case study, we use UW campus load data and solar data from the National Renewable Energy Laboratory (NREL)~\cite{NREL_solar_data} to build the known distribution $\mathbb{P}_Z$. Fig.~\ref{fig:load_scenarios_week} showcases twenty weekly load profiles generated via the GANs-based scenario generation technique.

\subsection{Interpreting the NWAs planning problem solution}

Our NWAs planning problem produces recommendations of optimal NWA capacity $\bd \phi_i^*$ and the timing of capacity expansion $\delta^*$. Naturally, the planner should re-evaluate the plan before implementation if conditions change. For instance, if the solution to the planning problem suggests that capacity should be expanded in 10 years, the planner should re-evaluate the plan before implementation, e.g., in 8 or 9 years. Similarly, if the plan suggests installing a NWA sometime in the future, the planner should re-evaluate using new data and conditions before implementation.   

Additionally, the proposed method produces optimal operating decisions $\bd x_i^*$. Their main purpose is to account for NWAs operational cost and benefits as these are crucial to determine optimal NWA investments. The operating decisions $\bd x_i^*$ are not meant to be implemented during real-time operation. Instead, the NWAs real-time operating strategy should be formulated using short-term predictions (e.g., minutes to days) and specialized control algorithms.

\section{Solution techniques} \label{sec:DWDA}
The non-convexity and high-dimensionality of Problem~\eqref{eq:CVX_problem} presents computational challenges that existing solvers cannot directly handle. Consider that for a time step length of 1 hour and a planning horizon of 20 years, the dimensionality of the sets $\bd{\mathcal{X}}_i(\bd{\phi}_i)$ ranges from roughly $175,000$ for the simplest cases (e.g., solar PV or EE) to more than half a million for the more complex ES case. Considering all four NWAs and the robust formulation, the problem in~\eqref{eq:CVX_problem} has roughly $2,000,000$ variables and constraints. While it is possible to attempt to ``brute-force''~\eqref{eq:CVX_problem} by formulating it as a mixed-integer program, the dimensionality of the problem is likely to doom such endeavor.   

\subsection{Technique 1: Dantzig-Wolfe Decomposition for NWAs}

Instead, we decompose Problem~\eqref{eq:CVX_problem} into $|\bd{\mathcal{N}}|$ subproblems using the DWDA to handle the dimensionality issue. Each NWA falls into a subproblem while a low-dimensional master problem handles the demand charge and the present cost of capacity expansion.

In our case studies, every subproblem is tractable. However, subproblem tractability is not necessarily true for more complex NWAs. However, it is possible to take advantage of the structure of the subproblems to solve them. For example, the investment variables $\bd \phi_i$ couple operations sub-subproblems that encompass a sufficiently small time horizon (e.g., a year or a month). Bender's decomposition is a suitable technique for problems coupled by variables~\cite{conejo2006decomposition}.

The master problem inherits the non-convexities of~\eqref{eq:CVX_problem}. We decompose the master problem and find its solution by solving a small number of small-scale linear programs. The details on this decomposition technique are found in Appendix~\ref{app:DWDA}.

\subsection{Technique 2: sequential solving}

The standard implementation of the DWDA may exhibit slow convergence rates due to a phenomenon called the ``tailing-off effect''~\cite{lubbecke2010column}. Fast convergence of control and short-term planning problems is critical since the delivery of solutions is time-sensitive.  In contrast, slow convergence of a long-term planning problem is not a big issue since long-term plans are typically implemented months or even years later. Nevertheless, even in long-term planning contexts, reasonable amounts solve times are important. Therefore, we provide an alternative solution technique that may converge faster than the DWDA in some instances. It is worth noting that, given enough computational resources, both techniques converge to the same objective value.   

\begin{algorithm}
    
    \KwIn{ $\mbox{P}$, $|\bd{\mathcal A}|$ }
    \KwOut{  $\mbox{objective value}$ }
    $j \leftarrow0$ \\
    \While{$j\le|\bd{\mathcal A}|$ and  $\mathrm{objective \;value}=\emptyset$}{
        \If{$\mbox{P}(j)$ is feasible}{
            
            \If{$j>0$ and $\mbox{P}(j)>\mbox{P}(j-1)$}{

                $\mbox{objective value}\leftarrow\mbox{P}(j-1)$  \\

            }
            
            \ElseIf{ $j =|\bd {\mathcal A}|$ }{
                
                $\mbox{objective value}\leftarrow \mbox{P}(j)$  \\
                
            }
        }
        
        \Else{
            
            $\mbox{objective value}\leftarrow\mbox{P}(j-1)$  \\
            
        }
        $j\leftarrow j+1$
        
    }
    \caption{Algorithm to sequentially solve the NWAs planning problem. The function $\mathrm P(j)$ represents Problem~\eqref{eq:CVX_problem} with the variable $\delta = j$.  }
    \label{alg:MP_solve}
    
\end{algorithm}

Technique 2 consists on convexifying Problem~\eqref{eq:CVX_problem} by fixing the year of capacity expansion $\delta$. We sequentially solve Problem~\eqref{eq:CVX_problem} for $\delta \in \{0,\dots, | \bd{\mathcal{A}}| \}$. It is not necessary to solve for $\delta$ in every year of the planning horizon. Since Problem~\eqref{eq:CVX_problem}'s objective is convex, we can stop once we find that the solution with $\delta = x $ is larger than the solution with $\delta=x + 1$ or $\delta$ is so large that a feasible solution can not be found as detailed by Algorithm~\ref{alg:MP_solve}. While this alternative may converge faster, it is not as scalable as the DWDA approach.

\subsection{Practical considerations}

Practitioners and intended end-users may not have the appropriate training in optimization and mathematics to implement our method. However, in commercial deployment, the optimization and mathematical models would run in the background. Thus, the end user does not necessarily need to make ``under-the-hood'' changes to the mathematics behind our model. We envision a commercial implementation of our method to ask the user for major economic and technical parameters of the system to be planned without the need to know any kind of optimization or complex mathematics.

\section{Case study: non-wire alternatives for the University of Washington}
\label{sec:case_study}
In this section, we study a NWAs planning problem at the UW Seattle Campus, which expects to add 6 million square feet of new buildings (e.g., labs, classrooms, office space) during the next ten years~\cite{UWMP_2018}. The additional load from new buildings will likely require expanding the capacity to serve the campus.

 Seattle City Light (SCL) and the UW are considering several traditional solutions to manage the expected load increase. The traditional solutions include building a new feeder to campus or increasing the service voltage to sub-transmission levels. However, these solutions are hard to implement in Seattle's dense urban environment and come at an estimated cost in the order of $\$100$ million. Moreover, there is an increasing appetite by SCL, the Washington State government, and the UW to explore novel approaches such as NWAs. 
 
 In this case study, we assume that the UW bears the investment and operating costs of the NWAs and that the cost of capacity expansion is shared by the UW and SCL. Furthermore, we assume that both parties are interested in achieving the least-cost solution.

\subsection{Data} \label{sec:data}

\begin{table}[]
    \centering
    \caption{Non-wire alternatives parameters}
    \label{table:NWA_parameters}
    \begin{tabular}{@{}|l|l|l|@{}}
        \toprule
        \rowcolor[HTML]{C0C0C0} 
        \textbf{Parameter}                      & \textbf{Value}            & \textbf{Source}           \\ \midrule
        \rowcolor[HTML]{EFEFEF} 
        \multicolumn{3}{|l|}{\cellcolor[HTML]{EFEFEF} \textbf{Energy Efficiency}} \\ \midrule
        Investment cost function                & piece-wise linear                       & \cite{brown2008us}\footnote{Adjusted to 2017 dollars and according to the assumption that buildings in Seattle are more energy efficient than the national average (due to Seattle having some of the strictest energy codes in the nation).}                \\ \midrule
        \rowcolor[HTML]{EFEFEF} 
        \multicolumn{3}{|l|}{\cellcolor[HTML]{EFEFEF}\textbf{Demand response}}                          \\ \midrule
        Investment cost                      & $\$200$/kW                  &      \cite{piette2015costs}                     \\ \midrule
        Efficiency coefficient               & $1.1$                      & modeling assumption               \\ \midrule
        \rowcolor[HTML]{EFEFEF} 
        \multicolumn{3}{|l|}{\cellcolor[HTML]{EFEFEF}\textbf{Energy Storage}}                           \\ \midrule
        Investment cost                         & $\$250$/kWh                 & \cite{Irena_ES}                    \\ \midrule
        Ch./dis.  efficiency            & $0.97/0.95$                 & \cite{Sarker_Optimal_2017}                 \\ \midrule
        Degradation coefficient                 & $0.028$ kWh/kW              & \cite{Sarker_Optimal_2017}                  \\ \midrule
        Energy-to-power ratio                   & $4$                         & \cite{lazard}                     \\ \midrule
        \rowcolor[HTML]{EFEFEF} 
        \multicolumn{3}{|l|}{\cellcolor[HTML]{EFEFEF}\textbf{Solar photovoltaics}}                                 \\ \midrule
        Investment cost                         & $\$2,000$/kW                     & \cite{fu2016nrel}                  \\ \midrule
        Production profile                      & -                       & \cite{NREL_solar_data}               \\ \bottomrule
    \end{tabular}

\end{table}

\begin{figure} 
\centering
\includegraphics[width=0.5\textwidth]{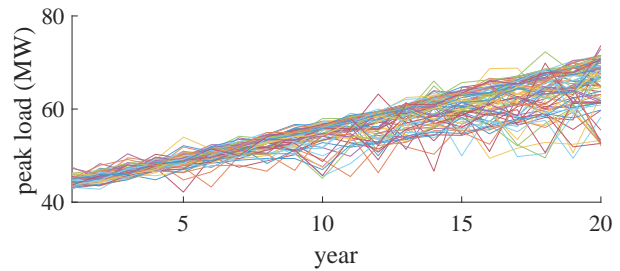}
\caption{Campus peak load scenarios. Note that some scenarios feature year-to-year peak load changes that are far from the average load growth. We include these low-likelihood scenarios to capture worst-case scenarios.  }\label{fig:load_data}
\end{figure}

 Table~\ref{table:NWA_parameters} summarizes the main parameters of the NWAs. We assume that, if installed, the NWAs begin operations at year 1. The NWAs value streams considered in this case study are energy arbitrage, peak-load reduction, and capacity expansion delay.   
 
 We assume that the cost of substation upgrades is $\$100$ million and adhere to a standard SCL planning horizon of 20 years~\cite{SCL_planning_practices_2018}.  We assume a yearly discount rate of $7\%$ and SCL tariffs for high-demand customers~\cite{SCL_tariff}. The N-1 pre-expansion capacity limit of the substation that serves the UW campus is 60 MW, i.e., if the substation loses one feeder, it can still serve 60 MW.

 We use NREL PV output data from a site near Seattle to generate PV scenarios~\cite{NREL_solar_data} and SCL campus load data from the years 2011 to 2016 to generate load scenarios.  Furthermore, we incorporate SCL's projected load growth of 1.5\% to 3.5\%  to the scenario-generation algorithm. Fig.~\ref{fig:load_data} shows the yearly peak load of each scenario. 

\subsection{Long-term planning results} \label{subsec:long-term planning results}

The main ``knob'' available to tune the results of the NWAs planning problem is the uncertainty protection level $\Gamma$ (see Appendix~\ref{app:robust} for more details). In this case study, we vary a single protection parameter $\Gamma \in [0,1]$ that determines the protection level of all uncertain parameters: load, PV generation, and EE reduction.

We interpret $\Gamma$ as follows. Suppose that we expect the load at a point in time to be within $50\pm5$ MW. Then, with a protection level of $\Gamma$, we optimize for the worst-case realization in the range $50\pm5\cdot\Gamma$ MW. For instance, with a protection level of $0.5$, the optimization problem considers load realizations within $50\pm2.5$ MW. With higher $\Gamma$'s, the optimization problem considers a broader range of possibilities and thus produce more robust solutions. 

\begin{figure} 
    \centering
    \includegraphics[width=0.45\textwidth]{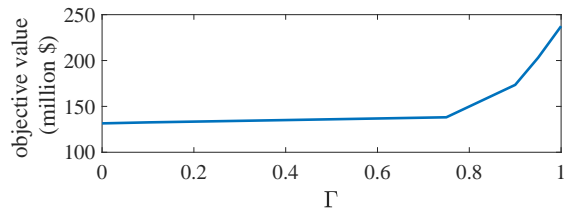}
    \caption{NWAs planning problem objective value as a function of $\Gamma$. } 
    \label{fig:planning_problem_objective}
\end{figure}

\begin{figure} 
    \centering
    \includegraphics[width=0.5\textwidth]{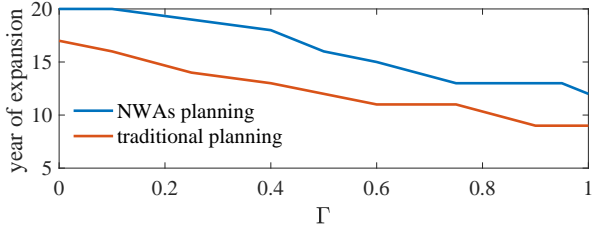}
    \caption{ Year of expansion for traditional and NWAs-base planning as a function of the protection level $\Gamma$. Postponing a $\$100$ million investment from year 9 to year 14 represents savings (with $\Gamma = 0.9$) of close to $\$13$ million.}
    \label{fig:YOE}
\end{figure}

However, more robust solutions represent higher costs. As shown in Fig.~\ref{fig:planning_problem_objective}, higher values of $\Gamma$ produce more expensive solutions to the NWAs planning problem. The solutions are more expensive in part because, as shown in Fig.~\ref{fig:YOE}, expanding capacity earlier (at a higher present cost) produces more robust solutions.

\begin{figure} 
    \centering
    \includegraphics[width=0.45\textwidth]{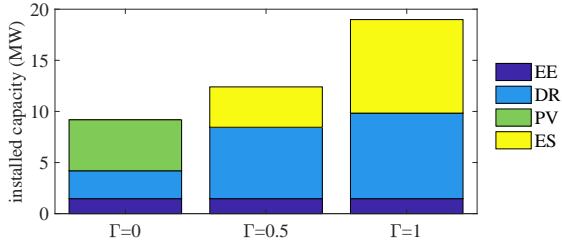}
    \caption{ Installed NWAs  capacities for three values of $\Gamma$: $0$, $0.5$, and $1$. In the case study, we assume that the NWAs are installed and start operation at year 1 of the planning horizon. That is, the year of installation is not a result of the optimization problem. } 
    \label{fig:NWA_Capacities}
\end{figure}

The level of protection also impacts the optimal mix of NWAs. As shown in Fig.~\ref{fig:NWA_Capacities}, the solutions for $\Gamma= 0.5$ and $\Gamma= 1$ favor technologies that, per our assumptions, do not face uncertainty (i.e., DR and ES). Conversely, the solution for $\Gamma=0$ ignores risks associated with uncertain PV production and thus favors the installation of solar PV. Higher levels of protection can be interpreted as the optimization assigning lower dependability to PV generation. Conversely, lower protection levels can be interpreted as the optimization algorithm being more ``optimistic'' about the dependability of PV generation. 

There are three main factors that affect the amount of PV installed capacity suggested by the solution to the NWAs planning problem. The first one, and as suggested previously, is the protection level. A solution with low protection level favor technologies that are subjected to uncertainty (in the case of PV, solar radiation levels). The second one is the correlation coincidence of load and solar production patterns. In a NWAs context, a higher coincidence translates into a greater ability of PV to offset peak load and thus mitigate the need for capacity expansion. The third and perhaps more obvious factor is the PV cost. Everything else equal, lower costs make PV investments more attractive.

Notice that as shown in Fig.~\ref{fig:NWA_Capacities}, the planning results suggest installing both DR and ES for protection levels $\Gamma = 0.5$ and $\Gamma = 1$. While ES and DR are to some extent substitute technologies, they are not perfect substitutes. In our model, ES can store charge as long as needed. DR, on the other hand, is not as flexible - it is limited to shifting load to an adjacent time period.  While the technical characteristics favor ES over DR, DR installations are cheaper in our model: 1 kW of DR-enabled load costs \$200 while a 1 kW/4 kWh battery costs \$1000. Roughly speaking, the optimization algorithm selects DR to shift load to adjacent time periods (since it is cheaper) and the more expensive ES to shift load across larger periods of time.   

\subsection{NWAs solution assessment via Monte Carlo simulations}

We perform Monte Carlo simulations on the possible realization of load and solar scenarios to asses the performance of the NWAs planning solution. 
\begin{figure}
    \centering
    \includegraphics[width=0.5\textwidth]{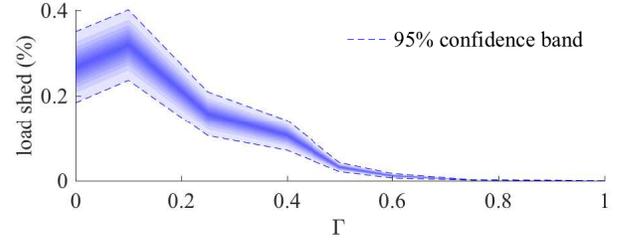}
    \caption{ Load shed (as percentage of total load) as a function of $\Gamma$.}
    \label{fig:load_shed}
\end{figure}

Perhaps the most dreaded consequence of planning a system without enough spare capacity, i.e., a non-robust system, is load-shedding. Fig.~\ref{fig:load_shed} shows the probability density of shedding load as a function of the protection level of the NWAs planning problem solution.  For the most part and as expected, load shedding decreases with $\Gamma$. Note that load-shedding (also known as energy-not-served) is only one metric of reliability. Other metrics such as loss-of-load probability or the customer average interruption duration index can easily be incorporated in the proposed framework. 

Load-shedding at U.S. university campuses is typically deemed undesirable. Thus, the UW campus should plan for NWAs with a $\Gamma$ guarantees no load shedding. That is, the planning problem should be solved with $\Gamma \approx 1$.

However, not every load requires such a high level of reliability.  The maximum price a load would be willing to pay not to be disconnected is known as the value of lost load (VOLL). For instance, the study in~\cite{LEAHY20111514} estimates that residential loads in Northern Ireland are willing to pay up to $18\euro$/kWh\footnote{In 2007 $\euro$.} to avoid load disconnections.

\subsection{How to select the level of protection?}

Although the protection level $\Gamma$ can be intuitively interpreted as in Section~\ref{subsec:long-term planning results}, we believe $\Gamma$ is too abstract for an appropriate value to be easily determined by most practitioners. Thus, we propose a methodology to maps a planner's VOLL (a concrete economic parameter) to the more abstract $\Gamma$. For a given VOLL, we select the value of $\Gamma$ that minimizes the sum of investment costs, expected energy costs, expected peak demand charges, and the expected cost of lost load. 

\begin{figure} 
    \centering
    \includegraphics[width=0.5\textwidth]{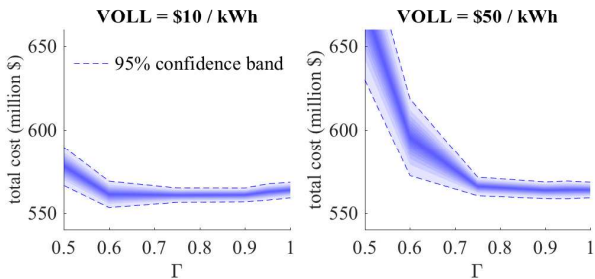}
    \caption{ Total as a function of $\Gamma$ for two different values of lost load.}
    \label{fig:cost_VOLL}
\end{figure}

\begin{figure} 
    \centering
    \includegraphics[width=0.45\textwidth]{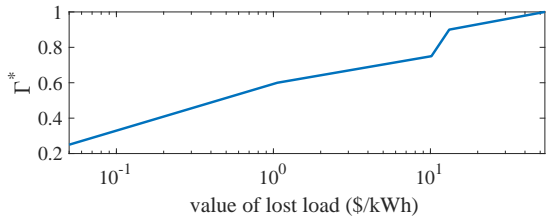}
    \caption{ Optimal protection level as a function of VOLL. We define the optimal protection level $\Gamma^*$ for a given VOLL as the $\Gamma$ that minimizes the expected total cost.  }
\end{figure}
For instance Fig.~\ref{fig:cost_VOLL} shows the total costs for two different VOLL. The left-hand plot shows the probability density of the total cost for $\mathrm{VOLL}=10$/kWh and its expected value reaches a minimum at $\Gamma\approx0.6$. The right-hand plot, on the other hand, shows the probability density of the total cost for $\mathrm{VOLL}=50$/kWh and its expected value reaches a minimum at a $\Gamma$ closer to $1$. Thus, a customer whose VOLL is $=10$/kWh would plan using $\Gamma\approx 0.6$ while a customer with a VOLL of $=10$/kWh would plan using $\Gamma\approx 1$.  We plot the value of $\Gamma$ that minimizes the expected total cost, $\Gamma^*$, as a function of VOLL in Fig.~\ref{fig:cost_VOLL}. Then, one can graphically map its VOLL to the appropriate $\Gamma$ to use in the NWAs planning problem.

One might ask, why not use the VOLL and minimize the expected total cost via stochastic optimization? That would be a good approach except for the fact that stochastic optimization is computationally more expensive than our robust optimization approach\footnote{The computational cost of a stochastic optimization problem increases with the number of scenarios. On the other hand, the size of robust optimization problems remains constant with the number of scenarios.}. Thus, the approach outlined in this paper is friendly to limited computing resources.

\subsection{Computational analysis}

This case study was carried out in a personal computer with 8 GB of installed memory (RAM) and running on an Intel\textregistered \: Core\texttrademark \:  i5-82500 CPU @ 1.60GHz. We implemented the optimization problems in Julia~\cite{Julia} and solved them using Gurobi 8.1.0.

\begin{figure} 
    \centering
    \includegraphics[width=0.45\textwidth]{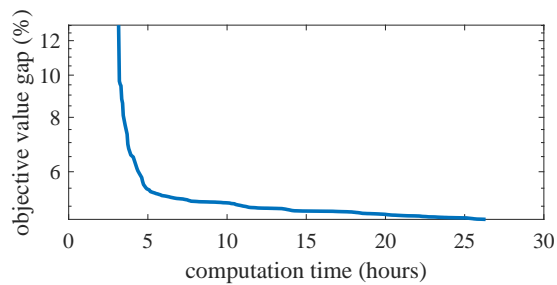}
    \caption{Gap between objective value and the optimal value as a function of computation time with technique 1 (DWDA-based) and $\Gamma=1$. }
    \label{fig: objective_vs_time}
\end{figure}

Fig.~\ref{fig: objective_vs_time} shows the percentage gap between the objective value of the NWAs problem and its optimal value obtained through the proposed DWDA-based technique. After 26 hours, the best solution (given by the last solution to the master problem) delivers a solution that is 4.66\% higher than the true optimal value. Note that the algorithm exhibits slow convergence, a well-documented issue known as the ``tailing-off'' effect. This issue can be addressed by tackling the degeneracy of the master problem as proposed by~\cite{desrosiers2011improved}.  Improving the convergence of the DWDA is beyond the scope of this work. 

\begin{figure} 
    \centering
    \includegraphics[width=0.45\textwidth]{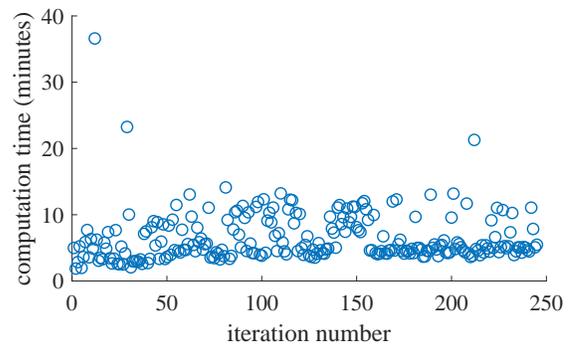}
    \caption{Computation time of each iteration as a function of iteration number with technique 1 (DWDA-based) and $\Gamma=1$. Each iteration takes an average of 6.44 minutes with a standard deviation of 3.71 minutes. }
    \label{fig: time_vs_iter}
\end{figure}

As shown in Fig.~\ref{fig: time_vs_iter}, the computation time of each iteration stays, on average, constant as the algorithm iterates. There are a couple of reasons, the first one is that the overwhelming majority (86.5\%) of the computation time is spent on the subproblems, whose size is constant throughout the algorithm. Only 13.5\% of the computation time is devoted to solving the master problem, whose size increases with the number of iterations. The second reason is that, in other to eliminate unnecessary information from the master problem, reduce its size and solve time, we remove proposals which have been assigned a weight of zero for at least 40 contiguous iterations. Removing zero-weight proposals from the master problem was suggested by Vanderbeck and Savelsbergh in~\cite{vanderbeck2006generic}.

\begin{figure} 
    \centering
    \includegraphics[width=0.45\textwidth]{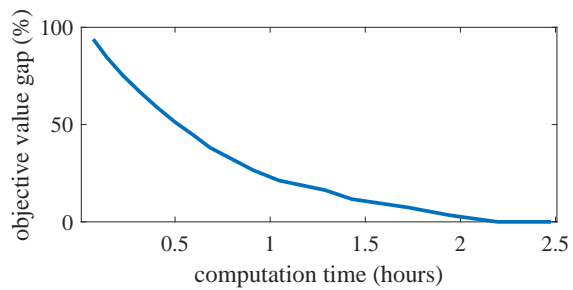}
    \caption{Gap between objective value and the optimal value as a function of computation time with technique 2 (sequential) and $\Gamma=1$. } \label{fig:obj_vs_time_seq}
\end{figure}

        Alternatively, the NWAs problem can be solved using the proposed technique 2 (sequential solving). As shown in Fig.~\ref{fig:obj_vs_time_seq}, the optimal solution to the problem is obtained in less than 2.5 hours. It is worth noting that while in this case technique 2 outperforms technique 1, it is less scalable than the DWDA-based technique. Thus, under certain circumstances (e.g., if considering a large number of NWAs and/or with insufficient computational resources) implementing technique 2 may not be possible.

\section{Conclusion and Future Work}
\label{sec:conclusion}
In long-term planning contexts, distributed energy resources (DERs) can serve as non-wire alternatives (NWAs) to traditional capacity expansion projects (e.g., increasing the capacity of a substation). In this paper, we formulate a NWAs planning problem that explicitly incorporates the time-value of delaying or avoiding traditional capacity expansion. Accounting for this additional value stream is crucial as it can define whether a DER project is economically viable or not.

We formulate the NWAs planning problem as a large-scale non-convex robust optimization problem. We tackle the size of the problem by proposing two techniques, both of which have pros and cons. The Dantzig Wolfe Decomposition Algorithm-based technique is scalable but can exhibit slow convergence rates. The second technique can provide faster convergence but has limited scalability. In practice, a utility can choose their preferred solution method according to the problem at hand.  Additionally, we present a case study that considers solar photovoltaic generation, energy efficiency, energy storage, and demand response as alternatives to substation/feeder upgrades at the University of Washington Seattle Campus.

	There are several open questions and research directions that we think are worth pursuing. For instance, we assume a radial system congested at a single point and disregard the downstream network. An interesting research direction is to consider a more general network (e.g., meshed or radial with multiple bottlenecks). Additionally, considering non-capacity NWAs such as network reconfiguration techniques is an interesting research direction. In this work, we focused on the capacity problem and ignored important phenomena such as power quality and transmission/distribution losses. More work should be done to integrate these phenomena into the NWAs planning problem. Regarding DER modeling, we believe that an interesting research question is: how to incorporate real-time control strategies and their performance into long-term planning problems? Regarding computation, we believe that research on how to increase the convergence and scalability of decomposed large-scale problems would be valuable.

\appendices

\section{Proof of Lemma~\ref{thm:IW_PV_CVX}} \label{sec:app_proof_lemma}
Stating that the function $\tilde{I}(\boldsymbol{l}^\mathrm{p})$ can be equivalently written as in \eqref{eq:IW_PV} and as the optimization problem in~\eqref{eq:IW_PV_CVX} is equivalent to stating that 
\begin{equation} 
\delta^* = \mathrm{CapEx}(\boldsymbol{l}^\mathrm{p}) \label{eq:equivalence}
\end{equation}
 where $\delta^*$ is the solution to Problem~\eqref{eq:IW_PV_CVX} and  $\mathrm{CapEx}(\boldsymbol{l}^\mathrm{p})$ is defined by Eq.~\eqref{eq:subs_upgrade_prob}. 

We prove Lemma by showing that Eq.~\eqref{eq:equivalence} in in fact holds. We start by noting that $\delta^* = \mathrm{CapEx}(\boldsymbol{l}^\mathrm{p})$ is a feasible solution to Problem~\eqref{eq:IW_PV_CVX}.  Suppose that $\epsilon$ is a small positive number. Then, the solution to Problem~\eqref{eq:IW_PV_CVX} cannot be $\delta = \mathrm{CapEx}(\boldsymbol{l}^\mathrm{p}) + \epsilon$ because it violates Constraint~\eqref{eq:NCVX_Const}.  Since we assume that $\rho$ is positive, the objective value given by $\delta = \mathrm{CapEx}(\boldsymbol{l}^\mathrm{p}) - \epsilon$ is larger than the objective value given by $\delta^*$. Since $\delta = \mathrm{CapEx}(\boldsymbol{l}^\mathrm{p}) + \epsilon$ is infeasible and $\delta = \mathrm{CapEx}(\boldsymbol{l}^\mathrm{p}) - \epsilon$ is suboptimal, the equality in~\eqref{eq:equivalence} and Lemma~\ref{thm:IW_PV_CVX} both hold. 
\QEDB
\section{Proof of Theorem~\ref{thm:NWAPP_cvx}}
\label{sec:appendix_proof_thm}
Using~\eqref{eq:IW_PV_CVX},~\eqref{prob:NWA_planning} can be written as 
\begin{equation}
\min_{\substack{\bd \phi_i \in \bd \Phi_i \\ \bd x_i \in \bd{\mathcal{X}}_i (\bd \phi_i) }}  \biggl\{ y + \min_{\delta \in \Delta } \frac{I^\mathrm{W}}{(1+\rho)^\delta}\biggr\}\label{eq:thm3_0}
\end{equation}
where $y=\sum_{i \in \bd{\mathcal{N}}} \left[C_i^\mathrm{O}(\bd x_i) + I_i^\mathrm{NW}(\bd \phi_i) \right]  +C^\mathrm{D}(\boldsymbol{l}^\mathrm{p}(\bd{x})) $ and $\Delta$ denotes the set defined by the constraints in~\eqref{eq:IW_PV_CVX}. Problem~\eqref{eq:thm3_0} is equivalent to 
\begin{equation}
\min_{\substack{\bd \phi_i \in \bd \Phi_i \\ \bd x_i \in \bd{\mathcal{X}}_i (\bd \phi_i) }} \min_{\delta\in\Delta } \biggl\{y +  \frac{I}{(1+\rho)^\delta}\biggr\} \label{eq:thm3_1}.
\end{equation}
Problem~\eqref{eq:thm3_1} is a nested optimization problem that whose inner variable is $\delta$ and its outer variables are $\bd x_i$ and $\bd \phi_i$. As shown in~\cite{boyd2004convex} the inner and outer variables can be minimized simultaneously in a single $\min$ function. Thus~\eqref{prob:NWA_planning} is equivalent to~\eqref{eq:CVX_problem} where~\eqref{eq:CVX_problem_c3} implement the functions $l^\mathrm{p}_a(\boldsymbol{x}) = \max_{t\in\bd{\mathcal{T}}}\{l_{a,t}(\boldsymbol{x})\}$ via half-planes. 
\QEDB

\section{Robust formulation of the NWAs planning problem} \label{app:robust}
We write \eqref{eq:CVX_problem}  in compact form as
\begin{subequations}
	\begin{align}
	\min_x \;&c^\top x \\
	\mathrm{s.t.}\; &Ax\le b
	\end{align} \label{eq:det_prob}
\end{subequations}
and assume that the uncertain parameters are contained in $A$. We implement uncertainty in $b$ by replacing the inequality with $\tilde A \tilde x\le 0$ where $\tilde A= [A,-b]$, $\tilde x = [x , \pi]^\top$, and $\pi=1$.

Let $a_{i,j}$ denote the $i,j$ element of $A$. Suppose that $a_{i,j}$ is a randomly distributed parameter with an unknown distribution that takes on values in $[\bar a_{ij} - \hat a_{ij}, \bar a_{ij} + \hat a_{ij}]$. Then, the robust counterpart of~\eqref{eq:det_prob} is
\begin{subequations}
	\begin{align}
	\min_{\substack{x,\; p_{i,j}\\ \; z_i,\; y_j}} \;& c^\top x \\
	\mathrm{s.t.}\; &\sum_j a_{i,j}x_j +z_i\cdot \Gamma_i + p_{i,j} \le b_i \; \forall i\\
	& z_i + p_{i,j}  \ge \hat{a}_{i,j}\cdot y_j \; \forall i , \;j\in J_{i} \\
	& -y_j \le x_j \le y_j\; \forall j\\
	& p_{i,j} \ge 0 \; \forall i, \; j\in J_i \\
	&y_j \ge 0 \; \forall j \\
	& z_i \ge 0 \; \forall i.
	\end{align} \label{eq:robust_det_prob}
\end{subequations}
Here, the set $J_i=\{j| \hat a_{i,j}>0\}$ and $\Gamma_i$ is allowed to take on values in $[0, |J_i|]$. The parameter $\Gamma_i$ adjusts the robustness of the solution and is known as the protection level of the $i^\mathrm{th}$ constraint ($\Gamma_i=|J_i|$ produces the most robust solution). The interested reader is encouraged to consult~\cite{bertsimas2003robust} for a more detailed explanation of~\eqref{eq:det_prob} and its robust counterpart~\eqref{eq:robust_det_prob}.

\section{Dantzig-Wolfe Decomposition for the NWAs planning problem} \label{app:DWDA}
The NWA subproblems are given by
\begin{equation*}
\min_{\substack{\bd \phi_i \in \bd \Phi_i \\
		\bd x_i \in \bd{\mathcal{X}}_i (\phi_i)}}  \;C_i^\mathrm{O}(\bd x_i) + I_i^\mathrm{NW}(\bd \phi_i) +\underbrace{ \sum_{a\in\bd {\mathcal A}} \sum_{t\in\bd{\mathcal T}} \pi^1_{a,t}\cdot l^i_{a,t}(\bd x_i)}_{\text{penalty term}}
\end{equation*}
for all $i\in\bd{\mathcal{N}}$. The objective of subproblem $i$ is composed of the operation cost $C_i^\mathrm{O}(\bd x_i)$, investment cost $I_i^\mathrm{NW}(\bd \phi_i)$, and a term that penalizes the load $l^i_{a,t}(\bd x_i)$ by $\pi^1_{a,t}$. The penalty coefficients $ \pi^1_{a,t}$ are the dual variables of the coupling constraints~\eqref{eq:master_problem_coupling} in the master problem.  The operating and investment decisions must be in their respective set of feasible solutions.

We write the master problem as
\begin{subequations}
	\begin{align}
	\min_{ \substack{\lambda_k,\; \delta \\ \boldsymbol{l}^\mathrm{p}=\{l_a^\mathrm{p}\}_{a\in\bd {\mathcal A}}} }& \left\{ \sum_{k=1}^K \lambda_k\cdot C^\mathrm{prop}_{(k)}+C^\mathrm{D}(\boldsymbol{l}^\mathrm{p})\!+\! \frac{I}{(1+\rho)^{\delta}} \right\}\\
	\mbox{s.t. } &  l^\mathrm{b}_{a,t} + \sum_{k=1}^K\lambda_{k} \cdot  l_{a,t,(k)}^{\mathrm{prop}}  \le l_a^\mathrm{p}\;\;\; (\pi^1_{a,t})  \;\forall \; a \in \bd {\mathcal A},\; t\in\bd{\mathcal{T} }\label{eq:master_problem_coupling}\\
	& l_{a}^\mathrm{p} \le \overline{l}\; \forall \; a<\delta \label{eq:master_problem_limit}\\
	& 0 \le \delta \le |\bd {\mathcal A}| \label{eq:master_problem_delta_lim}\\
	& \sum_{k=1}^K\lambda_{k} = 1 \;\;\; (\pi^2) \label{eq:master_problem_lambda_cvx} \\
	& \lambda_k \ge 0 \; \forall\; k=1,\dots,K \label{eq:master_problem_lambda_pos}
	\end{align} \label{eq:master_problem}
\end{subequations}
and its objective is to minimize the sum of three terms: a convex combination of $K$ cost proposals, peak demand charges, and the present cost of capacity expansion. The $k^\mathrm{th}$ cost proposal is
\begin{equation*}
C^\mathrm{prop}_{(k)}= \sum_{i \in \bd{\mathcal{N}}} C_i^\mathrm{O}(\bd x_{i,(k)}) + I_i^\mathrm{NW}(\bd \phi_{i,(k)})
\end{equation*}
where $x_{i,(k)}$ and $\phi_{i,(k)}$ represent optimal operating and investment decisions, respectively, for the $k^\mathrm{th}$ iteration. The positive variables $\lambda_k$ are the weights of each cost proposal. The coupling constraints are~\eqref{eq:master_problem_coupling}. We define the load proposals as
\begin{equation*}
l_{a,t,(k)}^{\mathrm{prop}}   = \sum_{i\in\bd{\mathcal{N}}} l^i_{a,t}(x_{i,(k)}).
\end{equation*}
Constraints~\eqref{eq:master_problem_limit} and ~\eqref{eq:master_problem_delta_lim}  originate from~\eqref{eq:CVX_problem_c4} and ~\eqref{eq:CVX_problem_c4.1}, respectively. Finally,~\eqref{eq:master_problem_lambda_cvx}  and~\eqref{eq:master_problem_lambda_pos} ensure that the sum of all $\lambda_k$'s equals one and that they are all non-negative. We skip the detailed description of the well-known DWDA. The interested reader is referred to~\cite{conejo2006decomposition} for an in-depth description and an implementation of the DWDA.

When decomposing~\eqref{eq:CVX_problem}, the non-convex Constraint~\eqref{eq:CVX_problem_c4} lies in master problem. We solve the master problem by solving $|\bd{\mathcal{A}}|+1$ linear problems. Let $\mathrm{P}(j)$ represent a function that sets $\delta=j$ in~\eqref{eq:master_problem} and solves for $\bd l^\mathrm{p}$ and $\lambda_k$. We solve $\mathrm{P}(j) \forall j \in \bd{\mathcal{A}}$ and select the smallest value of $\mathrm P$.  A concrete interpretation of $\mathrm{P}(j)$ is that capacity expansion happens at year $j$ and the peak load limit $\overline l$ is enforced from year $1$ through $j$. Note that the function $\mathrm{P}$ involves solving a small-scale linear problem. The number of variables in $\mathrm{P}$ is $K+|\bd {\mathcal A}|$ where $K$ is  in the order of a few hundred and $|\bd {\mathcal A}|$ is close to $20$.

\bibliographystyle{IEEEtran}
\bibliography{NWAsbiblio}

\end{document}